\documentclass[12pt,reqno]{amsart}
\usepackage{amscd,amsfonts,amssymb,amsmath,latexsym,array,hhline,amsthm}
\newtheorem{Theorem}{Theorem}
\newtheorem{Corollary}{Corollary}
\theoremstyle{remark}
\newtheorem{Remark}{Remark}

\begin{document}

\sloppy


\title[{Trigonometric Polynomials with Small Uniform Norm}]{Some Trigonometric Polynomials with Extremely Small Uniform Norm}

\author{Pavel G.~Grigoriev}
\address[]{P.G.~Grigoriev. Geogracom LLC}
\email{thepavel@mail.ru}

\author{Artyom O. Radomskii}
\address{A.O. Radomskii. Moscow Engineering Physics Institute}
\email[]{artrad@list.ru}

\date{June 8, 2014}


%
%
%
%
%
%
%

\maketitle

\begin{abstract}
{\small
An example of trigonometric polynomials with extremely small uniform norm is given.
This example demonstrates the potential limits for extension of Sidon's inequality
for lacunary polynomials in a certain direction.
}
\end{abstract}


We prove the following
\begin{Theorem}\label{Super Great Theorem 1}
Let $q>1$ and a sequence of naturals
$\{m_j\}_{j=1}^{\infty}$ satisfy
$m_{j+1}/m_{j}\geq q$ for all $j$.
Let another sequence of naturals
$\{d_j\}_{j=1}^{\infty}$ satisfy
$1\le d_j\le {m_{j+1}-m_j}$. Then there exists
a sequence of trigonometric polynomials
$\{\delta_{j}\}_{j=1}^\infty$
such that
\begin{gather}
\label{theoremineq0}
\delta_{j}(x)=\sum_{m_j\le s < m_j+d_j} c_s e^{isx},\\
\label{theoremineq1}
\frac18 \le \|\delta_j\|_1 \le \|\delta_j\|_\infty\le 7,
\\
\label{theoremineq2}
\Big\|\sum_{j=1}^N \delta_{j} \Big\|_\infty
\le \alpha
 +\beta \sqrt N +
\gamma \max_{1\le j\le N}\log_q \max\Big( \frac{m_j}{d_j}, \frac 1{\ln q},1\Big)
\end{gather}
for all $N=1,2,\dots$
with some positive absolute constants $\alpha$, $\beta$ and $\gamma$.
\end{Theorem}

This result improves the examples constructed by Grigoriev~\cite{ja1997}
and Radomskii~\cite{radom2013}, where roughly speaking Theorem~\ref{Super Great Theorem 1}
was proved for the case $m_j=2^j$, $d_j=[2^{j-j^\varepsilon}]$
and some small limitations.

\begin{Remark}
In Theorem~\ref{Super Great Theorem 1}
the example is constructed with the constants
$\alpha=316$,
$\beta=7\sqrt{2c_H}$,
$\gamma=210$,
where $c_H$ is the constant from the Carleson-Hunt inequality,
see~(\ref{major}) below.
Instead of the Carleson-Hunt result one could use
a well-known simpler inequality
(see~\cite{zygmund}, Ch.~13, Th.~1.17)
$$
\Big\|
\sup_{1\le j<\infty}
\big|
\sum_{k=1}^{m_j}
b_k e^{i kx}
\big|
\Big\|_2^2
\le A_q \sum_{k=1}^{\infty}
|b_k|^2
$$
whenever $m_{j+1}/m_{j}\geq q$
with a constant $A_q$ depending on $q$.
The Carleson-Hunt inequality is used here because we
would like to make the constants $\alpha$, $\beta$ and $\gamma$ independent of $q$.
Otherwise we do not try to optimize the constants in Theorem~\ref{Super Great Theorem 1},
e.g. for the case $q=2$  it is not difficult to repeat the arguments of our proof
and get some essentially better constants.
\end{Remark}

\begin{Remark}
{ In this paper we assume that the norms of $L_p(0,2\pi)$ are normalized, i.e.
$\|f\|_p=\big(\frac1{2\pi}\int_0^{2\pi} |f|^p \,d\mu\big)^{1/p}$ for $1\le p <\infty$,
where $\mu$ is the standard Lebesgue measure so that $\|1\|_p=1$.}
\end{Remark}

\begin{Remark}
One can deduce as a corollary a version of Theorem~\ref{Super Great Theorem 1}
with {\it real} polynomials $\delta_{j}$ with slight changes in the constants.
(We decided to prove the result for $\delta_{j}$ with {\it positive frequencies} and {\it real coefficients}.)
Also with minor changes in the proof one can prove a version of Theorem~\ref{Super Great Theorem 1}
with $\delta_{j}(x)=p_j(x)\cos m_j x$ with some real trigonometric polynomials $p_j$ such that $\deg p_j\le d_j$.
\end{Remark}

Taking Remark~2 into account we derive from
Theorem~\ref{Super Great Theorem 1} the following

\begin{Corollary}
 Let $0 \leq \varepsilon <1.$
Then there exists a sequence of real
trigonometric polynomials
$\{p_{j}\}_{j=1}^\infty$ such that
$\deg p_j\le  2^{j-j^\varepsilon}$,
$c_-\le\|p_j\|_1 \le \|p_+\|_\infty\le c_+$
and
$$
\Big\|\sum_{j=1}^{N} p_{j}(x)\cos 2^{j} x \Big\|_{\infty}\leq C N^{\max({\varepsilon},\frac12)}
$$
for all $N=1,2,\dots$ with some positive absolute constants $c_-$, $c_+$, $C$.
\end{Corollary}

The example constructed in Theorem~\ref{Super Great Theorem 1}
gives some  limits for the attempts to extend
the well-known property of lacunary polynomials
(Sidon's inequality)
$$
\Big\|\sum_{j=1}^N a_j \cos m_jx\Big\|_\infty\ge
c(q) \sum_{j=1}^N |a_j|
$$
with some $c(q)>0$
whenever $m_{j+1}/m_{j}\ge q>1$.
The question of possible extension of Sidon's inequality
with substituting $a_j \cos m_jx$
by $p_j(x) \cos m_jx$ with $p_j$ being
trigonometric polynomials of possibly large degree
was raised by Kashin and Temlyakov~\cite{kashtemlyak20072008}
and Radomskii~\cite{radom2011}
 in association with
 estimating the entropy numbers of certain functional classes.

 The scheme of proof follows that from Grigoriev~\cite{ja1997}, \cite{jadis}
and Radomskii~\cite{radom2013} with minor changes.
This technique invented in \cite{ja1997} could be called the
{\it pseudo stopping time method} since the idea was borrowed
from stochastic analysis
(see \cite{jadis} for explanations).

\textsc{Proof of Theorem~\ref{Super Great Theorem 1}.}
It is easy to notice that without loss of generality
we can consider only the case when
\begin{equation}\label{newlimitsDn}
\frac{m_j}{d_j}\ge\max\Big(1,\frac1{\ln q}\Big).
\end{equation}

By the famous result of Carleson and Hunt~\cite{Hunt}
the $L_2$-norm of the majorant of a trigonometric
sum
can be
estimated as
\begin{equation}\label{major}
\Big\|
\max_{0\le k\le n}
\big|
\sum_{s=0}^k b_s e^{isx}
\big|
\Big\|_2^2
\le
c_H
\sum_{s=1}^n |b_s|^2
\end{equation}
with an absolute constant $c_H>0$.

We construct the polynomials $\delta_{k_j}$
by induction with the constants
\begin{equation}\label{theconstants}
\alpha:=316, \quad
\beta:=7\sqrt{2c_H}, \quad
\gamma:=210.
\end{equation}
Let $\delta_{1}(x):=\exp(i m_1x)$.
On each inductive step
in order to construct $\delta_{n}$
having $\delta_{1},\dots,\delta_{{n-1}}$
from the previous steps
we denote
$$
S_j:=\sum_{k=1}^j\delta_j
$$
and define
\begin{align}
\label{En}
E_n^j
&:=\big\{x\in[0,2\pi): |S_j(x)|> \beta \sqrt{n}
\big\};\\
\label{Bn}
B_n
&:=\bigcup_{j=1}^{n-1} E_n^j;\\
\label{tildBn}
\widetilde{B}_n
&:=\bigcup_{1\le j\le n-a_n} O_{2\pi \frac{(n-k)^2}{d_n}} (E_n^j),\\
\label{a_n}
\text{where\ }
a_n
&:=45+30\log_q\frac{m_n}{d_n}
\end{align}
and
 $O_\varepsilon(X)$
denotes the
$\varepsilon$-neighborhood of the set $X$ on the semi-interval $[0,2\pi)$
with the circle metric.
{\it Note that the value $n-a_n$ need not to be neither integer nor positive.
If $n-a_n<1$, the union in (\ref{tildBn}) is understood as an empty set. Similarly
below
the sums like $\sum_{1\le j\le n-a_n}$
are supposed equal to zero whenever $n-a_n<1$.
}

Set
\begin{align}
\label{Lambda}
\Lambda_n
&:=\Big\{l\in\{1,\dots,d_n\} : 2\pi\frac{l}{d_n}\notin \widetilde{B}_n
\Big\};\\
\label{delta}
\delta_n(x)
&:=\frac1{d_n} \exp\Big\{i\big(m_n+\big[\frac{d_n-1}2\big]\big)x\Big\}
\sum_{l\in \Lambda_n}
K_{\big[\frac{d_n-1}2\big]}\Big(x- 2\pi\frac{l}{d_n}\Big),
\end{align}
where $K_d(x)=\frac2{d+1}\big(\frac{\sin\frac12 (d+1)x}{2\sin\frac x2}\big)^2 $
are the Fej\'{e}r kernels.
Clearly so defined $\delta_n$ satisfies the frequency constraint imposed by (\ref{theoremineq0}).

In order to verify the induction hypothesis we need to show that
\begin{equation}\label{Lambdaineq}
|\Lambda_n|>\frac14  d_n.
\end{equation}

By the induction hypothesis
we have (\ref{theoremineq1})
for  $\delta_{1},\dots,\delta_{{n-1}}$.
Applying the Chebyshev inequality for the majorant
$S^*_{n-1}:=\max\limits_{1\le j \le n-1}\big|\sum\limits_{s=1}^j \delta_{s}(x)\big|$
(see (\ref{En}), (\ref{Bn}))
and then using (\ref{major}), (\ref{theoremineq1})
and (\ref{theconstants})
we get
\begin{equation}\label{muB}
\mu B_n\le
\frac{2\pi \|S^*_{n-1}\|_2^2}{\beta^2 n}
\le
\frac{2\pi c_H\|S_{n-1}\|_2^2}{\beta^2 n}
\le
\frac {2\pi }{\beta^2  n} \sum_{j=1}^{n-1}\|\delta_{k_j}\|_\infty^2
\le
\frac {2\pi }{\beta^2  n} (n-1)7^2
<\pi.
\end{equation}

Let us denote by
$\mbox{Conn}(X)$
the number of connected components
of a set $X\subset[0,2\pi)$ in the circle topology.
Note that $|S_j|^2$
is a real trigonometric polynomials
of degree not exceeding $2(m_j+2[\frac{d_j-1}2])$
and therefore the equation $|S_j(x)|^2=\beta^2 n$
has not more than $4(m_j+2[\frac{d_j-1}2])$
roots. Note that $E_n^j$
is a finite union of open intervals which endpoints are the roots of $|S_j(x)|^2=\beta^2 n$.
Taking into account (\ref{newlimitsDn})
we conclude
\begin{equation}
\label{connEn}
\mbox{Conn}(E_n^j)
\le \frac12 \big|\{x: |S_j(x)|^2=\beta^2 n\} \big|
\le 2 \Big(m_j+2 \big[\frac{d_j-1}2\big] \Big)
<4m_j.
\end{equation}
This implies
\begin{align*}
\mbox{Conn}(\widetilde{B}_n)
&\le
\sum_{1\le j\le n-a_n} \mbox{Conn}(E_n^j)
<
4\sum_{1\le j\le n-a_n} m_j
\le 4\sum_{1\le j\le n-a_n} m_n q^{j-n}\\
&<
 4m_n\sum_{s\ge a_n}  q^{-s}
 =4m_n q^{-\lceil a_n \rceil}\frac{q}{q-1},
\end{align*}
where $\lceil x \rceil$ denotes the ceiling integer part of $x$.
Taking into account that
$a_n\le\lceil a_n \rceil$, $\ln q\le q-1$
and using (\ref{a_n}) and (\ref{newlimitsDn})
we proceed as
\begin{align*}
\mbox{Conn}(\widetilde{B}_n)
&<
4m_n q^{- a_n }\frac{q}{q-1}
 =4m_n q^{-45}\big(\frac{m_n}{d_n}\big)^{-30}\frac{q}{q-1}
=
4 d_n \big(\frac{m_n}{d_n}\big)^{-29}\frac{q^{-44}}{q-1}\\
&\le
4 d_n \big(\frac{m_n}{d_n}\big)^{-29}\frac{q^{-44}}{\ln q}
\le
4 d_n \big(\frac{m_n}{d_n}\big)^{-28} {q^{-44}}
\le
4 d_n \max\Big( \big(\frac{m_n}{d_n}\big)^{-28}, q^{-44}\Big).
\end{align*}
If
 $q\ge e^{1/2}$,
 then
$ q^{-44}\le e^{-22}$.
If
 $q\le e^{1/2}$,
 then
$ \big(\frac{m_n}{d_n}\big)^{-28}\le \big(\frac{1}{\ln q}\big)^{-28}\le 2^{-28}
$.
So we get
\begin{equation}
\label{conntildeBn}
\mbox{Conn}(\widetilde{B}_n)
<
4 d_n \max\Big( 2^{-28},e^{-22}\Big)
<\frac{d_n}8.
\end{equation}

Aggregating
(\ref{muB})
and
(\ref{connEn}) (see also (\ref{tildBn}))
we get
\begin{align*}
\notag
\mu \widetilde{B}_n
&\le\mu B_n+
2\sum_{1\le j\le n-a_n}
2\pi\frac{(n-j)^2}{d_n}\mbox{Conn}(E_n^j)
\le
\pi+
16\pi \sum_{1\le j\le n-a_n}
\frac{(n-j)^2}{d_n} m_j\\
&
\le
\pi+
16\pi \sum_{1\le j\le n-a_n}
\frac{(n-j)^2}{d_n} m_n q^{j-n}
=
\pi+
16\pi \frac{m_n}{d_n} \sum_{a_n\le s\le n-1}
s^2 q^{-s}.
\end{align*}

To proceed with the estimate of $\mu \widetilde{B}_n$
we are going to use the inequality
$$
\sum_{s=a}^\infty s^2q^{-s}
\le\frac{2a^2q^{3-a}}{(q-1)^3}
\qquad
\text{for\ }
a=1,2,\dots.
$$
This inequality one  could easily
deduce from the following identity
$$
\sum_{s=a}^\infty s^2q^{-s}
=\frac{q^{3-a}}{(q-1)^3}
\Big(
a^2+q^{-1}(1+2a-2a^2)+q^{-2}(a-1)^2
\Big).
$$
for all $q>1$ and $a=0,1,\dots$, which is not so difficult to verify.

Now proceed with the estimate of $\mu \widetilde{B}_n$
as
$$
\mu \widetilde{B}_n
<
\pi+
32\pi \frac{m_n}{d_n} \frac{\lceil a_n \rceil^2q^{3-\lceil a_n \rceil}}{(q-1)^3}
<
\pi+
32\pi \frac{m_n}{d_n} \frac{\lceil a_n \rceil^2q^{3-\lceil a_n \rceil}}{\ln^3q}.
$$
The function $x^2q^{-x}$
is decreasing for $x\ge 2/\ln q$.
One easily checks that (\ref{newlimitsDn}) and (\ref{a_n})
imply $a_n\ge 2/\ln q$
and therefore we can substitute $\lceil a_n \rceil$ for $ a_n $ in the right-hand side
above. So recalling
(\ref{a_n}) and (\ref{newlimitsDn}) again we continue as
\begin{align*}
\mu \widetilde{B}_n
&<
\pi+
32\pi \frac{m_n}{d_n} \frac{ a_n ^2q^{3- a_n }}{\ln^3q}\\
&=
\pi+
32\pi \big(\frac{m_n}{d_n}\big)^{-29} \frac{ (45+30\log_q\frac{m_n}{d_n}) ^2q^{-42 }}{\ln^3q}\\
&
\le
\pi+
32\pi \big(\frac{m_n}{d_n}\big)^{-26} (45+30\log_q\frac{m_n}{d_n}) ^2q^{-42 }\\
&
\le
\pi+
64\pi \big(\frac{m_n}{d_n}\big)^{-26} \Big\{45^2+30^2\Big(\frac{\ln\frac{m_n}{d_n}}{\ln q}\Big) ^2\Big\}q^{-42 }\\
&
\le
\pi+
64\pi  \Big\{45^2\big(\frac{m_n}{d_n}\big)^{-26}+30^2\big(\frac{m_n}{d_n}\big)^{-22}\Big\}q^{-42 }.
\end{align*}
If $q\ge e^{1/2}$, then
$$
\Big\{45^2\big(\frac{m_n}{d_n}\big)^{-26}+30^2\big(\frac{m_n}{d_n}\big)^{-22}\Big\}q^{-42 }
\le
(45^2+30^2)e^{-21 }
<
(64^2+64^2)2^{-21 }=2^{-8}.
$$
If $q\le e^{1/2}$, then $\frac{m_n}{d_n}\ge1/\ln q\ge2$ and
$$
\Big\{45^2\big(\frac{m_n}{d_n}\big)^{-26}+30^2\big(\frac{m_n}{d_n}\big)^{-22}\Big\}q^{-42 }
<
45^2 2^{-26}+30^2 2^{-22}<2^{-11}.
$$
So we finally conclude
\begin{equation}
\label{mutildB}
\mu \widetilde{B}_n
< \pi
+64\pi 2^{-8}
= \frac54 \pi.
\end{equation}

Now we are ready to prove
(\ref{Lambdaineq}).
Clearly, the number of elements in $\Lambda_n^c$
does not exceed the number of the intervals of type
$\big(\frac{2\pi l-\pi}{d_n},\frac{2\pi l+\pi}{d_n}\big)$
($l=1,\dots,d_n$) which intersect
$\widetilde{B}_n$ (see (\ref{Lambda})).
Such intervals can be split into two groups:
 those containing an edge point of $\widetilde{B}_n$
and those being included in $\widetilde{B}_n$.
There are not more than
$2 \mbox{Conn}(\widetilde{B}_n)\le d_n/8$ of the intervals of the first type
(see (\ref{conntildeBn})).
Denoting by $V$ the number
of the intervals of the second type
and recalling (\ref{mutildB})
we get
$$
V\frac{2\pi} {d_n}
\equiv
V\mu\big(-\frac{\pi}{d_n},\frac{\pi}{d_n}\big)
\le \mu\widetilde{B}_n<\frac54\pi.
$$
So we have
$$
|\Lambda_n|\ge d_n-V-\frac {d_n}8 >d_n-\frac 58 d_n-\frac {1}8 d_n=\frac14 d_n.
$$
Thus we proved (\ref{Lambdaineq}).

Our next goal is to verify
(\ref{theoremineq1}) for $\delta_{n}$.
Let us recall some properties of the
Fej\'{e}r kernels. It is well-known that for
all $-\pi\le x\le \pi$ and $d=1,2\dots$
\begin{gather}
\label{fejerprop}
K_{d-1}(x)\ge0,\qquad
\|K_{d-1} \|_1=1/2,
\\
\label{fejerineq}
K_{d-1}(x)\le\min\big(\frac{d}2, \frac{\pi^2}{2dx^2}\big)
< 5
\min\big(d, \frac{1}{dx^2}\big).
\end{gather}

Denote  by
$\mbox{dist}(x,y)$ the standard distance on the circle  between $x,y\in[0,2\pi)$
and let $\widetilde{d}_n:=\big[\frac{d_n-1}2\big]+1$. Clearly, $1\le d_n/\widetilde{d}_n\le 2$.
Using  (\ref{delta}), (\ref{Lambdaineq}), (\ref{fejerprop}) and
(\ref{fejerineq}) we verify (\ref{theoremineq1}) as follows
$$
\|\delta_{n}\|_1
=\frac1{d_n}
\Big\|
\sum_{l\in\Lambda_n}
K_{\widetilde{d}_n-1
}\Big(x- 2\pi\frac{l}{d_n}\Big)
\Big\|_1
=\frac{|\Lambda_n|}{d_n} \big\|K_{\widetilde{d}_n-1}\big\|_1 \ge\frac18
$$
and
\begin{align*}
\|\delta_{n}\|_\infty
&=
\frac1{d_n}
\Big\|
\sum_{l\in\Lambda_n}
K_{\widetilde{d}_n-1}\Big(x- 2\pi\frac{l}{d_n}\Big)
\Big\|_\infty\\
&\le
\frac1{d_n}
\Big\|
\sum_{l=1}^{d_n}
5
\min\big(\widetilde{d}_n, \frac{1}{\widetilde{d}_n\,\mbox{dist}(x,2\pi\frac l{d_n})^2}\big)
\Big\|_\infty\\
&\le
5\Big(1+2
\sum_{s=1}^\infty\frac1{\widetilde{d}_n^2(2\pi\frac s{d_n})^2}
\Big)\\
&=5\Big(1+
\frac{d_n^2}{2\pi^2 \widetilde{d}_n^2}
\sum_{s=1}^\infty\frac1{s^2}
\Big)\\
&=5\Big(1+
\frac{d_n^2}{2\pi^2 \widetilde{d}_n^2}
\frac{\pi^2}6
\Big)<7.
\end{align*}

Now to complete the proof it remains to verify (\ref{theoremineq2})
with the constants
(\ref{theconstants}), i.e. to show that
$$
|S_n(x)|\le\alpha
 +\beta \sqrt n +
\gamma \max_{1\le j\le n}\log_q \max\Big( \frac{m_j}{d_j}, \frac 1{\ln q},1\Big)
\quad\text{for each\ }x\in[0,2\pi).
$$
Set
$$
\tau(x):=\max\{t=1,\dots,n-1:|S_t(x)|\le \beta\sqrt{ n}\}.
$$
If
$
\tau(x)\ge n-a_n
$,
then using (\ref{theoremineq1}) and (\ref{a_n}) we get
\begin{align*}
\notag
|S_n(x)|
&\le |S_{\tau(x)}(x)|
+\sum_{t=\tau(x)+1}^{n}
|\delta_{t}(x)|\le
\beta \sqrt{ n}+ 7a_n\\
&
=\beta \sqrt{ n}+ 7\Big(45+30\log_q\frac{m_n}{d_n}\Big)
=\alpha-1
+\beta \sqrt{ n}+ \gamma\log_q\frac{m_n}{d_n}.
\end{align*}
If
$
\tau(x)< n-a_n
$,
then
\begin{align}
\notag
|S_n(x)|
&\le |S_{\tau(x)}(x)|
+\sum_{t=\tau(x)+1}^{\tau(x)+a_n}
|\delta_{t}(x)|+
\sum_{t=\tau(x)+a_n+1}^n
|\delta_{t}(x)|\\
\label{tau3.22}
&\le
\alpha-1
+\beta \sqrt{ n}+ \gamma\log_q\frac{m_n}{d_n}
+\sum_{t=\tau(x)+a_n+1}^n
|\delta_{t}(x)|.
\end{align}

It remains to estimate the last term in (\ref{tau3.22}).
Since $\tau(x)+a_n+1\le t\le n$,
we have $x\in E_n^{\tau(x)+1}\subset E_t^{\tau(x)+1}$
(see (\ref{En})).
Therefore, by the definition of $\widetilde{B}_t$
(see (\ref{tildBn}))
we conclude
$$
\inf_{y\in [0,2\pi)\setminus \widetilde{B}_t}
\mbox{dist}(x,y)
\ge\frac{2\pi}{d_t}|t-\tau(x)-1|^2.
$$
Consequently, for each $l\in \Lambda_t$ we have (see (\ref{Lambda}))
$$
\mbox{dist}\Big(x,2\pi\frac{l}{d_t}\Big)
\ge\frac{2\pi}{d_t}|t-\tau(x)-1|^2.
$$
Using (\ref{fejerineq}) again and applying the trivial estimates
$\sum_{s=K+1}^\infty s^{-2} <K^{-1}$
and
$\sum_{s=K}^\infty s^{-2}< 2K^{-1}$
we conclude
\begin{align*}
|\delta_{t}(x)|
&
\le \frac1{d_t}\sum_{l\in\Lambda_t} K_{\widetilde{d}_t-1}\Big(x- 2\pi\frac{l}{d_n}\Big)
\le
\frac 5{d_t}\sum_{l\in\Lambda_t} \min\Big(\widetilde{d}_t,\frac1{\widetilde{d}_t\mbox{dist}(x,2\pi\frac{l}{d_t})^2}\Big)\\
&\le
10\sum_{s=|t-\tau(x)-1|^2}^\infty
\frac1{\widetilde{d}_t d_t \big(2\pi\frac s{d_t}\big)^2}
=
\frac {5d_t} {2\pi^2\widetilde{d}_t} \sum_{s=|t-\tau(x)-1|^2}^\infty
\frac{1}{   s^2}\\
&<\sum_{s=|t-\tau(x)-1|^2}^\infty
\frac{1}{   s^2}
<
\frac 2{|t-\tau(x)-1|^2}.
\end{align*}
Now we can estimate the last term in
(\ref{tau3.22}) as
\begin{equation*}
\sum_{t=\tau(x)+a_n+1}^n|\delta_{t}(x)|
<
\sum_{t=\tau(x)+a_n+1}^n \frac2{|t-\tau(x)-1|^2}
=
\sum_{s=a_n}^\infty \frac2{s^2}< \frac2 {a_n-1}<1.
\end{equation*}
Using the last estimate in (\ref{tau3.22}) we
get (\ref{theoremineq2}) for $S_n$.
This completes the proof.

\end{document}